\input amstex
\magnification=\magstep1
\documentstyle{amsppt}
\vsize 45pc

\NoBlackBoxes
\TagsAsMath

\define\enddemoo{$\square$ \enddemo}

\define\Hom#1#2#3{\text{Hom}_{#1}({#2},{#3})}
\define\End#1#2{\text{End}_{#1}({#2})}
\define\Ann#1#2{\text{Ann}_{#1}({#2})}

\NoRunningHeads

\document
\topmatter
\title Tight Closure Commutes with Localization in Binomial rings.
\endtitle
\author   Karen E. ~Smith $^{*}$ 
\endauthor
\address 
 Mathematics Department, University of Michigan, Ann Arbor, Michigan 48109
\endaddress 
\thanks {$^{*}$ Supported by the National Science Foundation and the Alfred P. Sloan
Foundation.}
\endthanks
\keywords tight closure, localization, binomial ideals, toric rings, semigroup rings
\endkeywords
\subjclass 13A35
\endsubjclass
\abstract It is proved that tight closure commutes with localization 
in any domain  which has a module finite extension in which tight closure is
known to commute with localization. It follows that tight closure commutes with localization
in binomial rings, in particular in semigroup or toric rings.
\endabstract
\endtopmatter

Tight closure, introduced in \cite{HH1}, is a closure operation performed on ideals in a
commutative, Noetherian ring containing a field. Although it has many important applications
to commutative algebra and related areas, 
some basic questions about tight closure remain open. 
 One of the biggest such open problems  is whether or not tight closure
 commutes with
localization. For an exposition of tight closure, including the
localization problem, see \cite{Hu}.

 The purpose of this paper is give  a simple proof that tight closure commutes with
localization for  a  class of rings that includes semi-group rings  and slightly more
generally, the so-called ''binomial" rings, that is, those that are a quotient of a
finitely generated algebra by an ideal generated by binomials.   Although this is a modest
class of rings,  the result here includes coordinate rings of all affine toric varieties,
and is general enough to   substantially improve previous  results of a number a
mathematicians, including  Moty Katzman, Will Traves, Keith Pardue and Karen Chandler, Irena
Swanson and myself. Thanks to Florian Enescu for pointing out 
a missing hypothesis in an earlier version.

\medskip
Throughout this note,  $R$ denotes a Noetherian commutative
 ring of prime characteristic.

We say that {\it tight closure commutes with localization in $R$\/}
if, for all ideals $I$ in $R$ and all multiplicative systems $U$ in $R$,
$I^*R[U^{-1}] = (IR[U^{-1}])^*$.

The main result is the following:

\proclaim{Theorem}
Let $R$ be a ring with the following property: for each minimal prime  $P$ of 
$R$, the quotient $R/P$ has a finite 
integral extension domain in which tight closure
commutes with localization. Then tight closure commutes with localization in $R$. 
\endproclaim

\proclaim{Lemma 1}
If tight closure commutes with localization in $R/P$ for each minimal
prime $P$ of $R$, then tight closure commutes with localization in $R$.
\endproclaim

\demo{Proof}(See also {\cite{AHH 3.8}})
It follows immediately from the definitions that $I^*R[U^{-1}] \subset (IR[U^{-1}])^*$
in general. To proof the reverse inclusion, 
let $\frac{z}{1}$ be in $(IR[U^{-1}])^*$. Then modulo each minimal prime of $R$ disjoint 
from $U$, the image of $\frac{z}{1}$ is in $(IR/P[U^{-1}])^* = (IR/P)^*R/P[U^{-1}]$.
For each minimal prime $P_i$ disjoint from $U$, there is some $u_i$ in $U$ such that 
$u_iz$ is in $(IR/P_i)^*$. By setting 
$u = \prod u_i$,
we see that $uz \in (IR/P)^*$ for all minimal primes $P$ of $R$ disjoint from $U$.
On the other hand, for each minimal prime $Q$ that meets $U$, pick any  any element
in $Q \cap U$ and let $v$ be their product.  Now $uvz \in (IR/P)^*$ for every minimal
prime $P$ of $R$, which means that $uvz \in I^*$ in $R$. But because $uuv \in U$, we
conclude
that $\frac{z}{1} \in I^*R[U^{-1}]$ and tight closure commutes with localization in $R$.
\enddemoo

\proclaim{Lemma 2}
Let $R$ be a domain  that has a finite extension domain $S$ in which tight closure commutes
with localization. Then tight closure commutes with localization in $R$.
\endproclaim

\demo{Proof}
Let $z\in R$ and let  $I \subset R$.   If $\frac{z}{1} \in (IR[U^{-1}])^*$, then 
the equations $c(\frac{z}{1})^q \in I^{[q]}R[U^{-1}]$ hold also after expansion to
$S[U^{-1}]$, so that $\frac{z}{1} $ is in $(IS[U^{-1}])^*$. Since tight closure commutes 
with localization in $S$, we have $\frac{z}{1}$ is in $(IS)^*S[U^{-1}]$. This means that 
there exists a $u \in U$ such that $uz \in (IS)^*$. But in general it is easy to see that,
$(IS)^*
\cap R
\subset I^*$ for any ideal $I$ of $R$ and any integral extension domain $S$ of $R$, so
$uz \in I^*$.
(Indeed, by tensoring with the fraction field of $R$ and then clearing denominators,
one can find an $R$-linear map $S @>\phi>> R$  which sends a fixed non-zero $c$  in $S$ to 
some non-zero element $d$ of $R$.  Then we can apply $\phi$ to the equations
expressing $c(zu)^q \in
I^{[q]}S$ to get equations expressing $d(uz)^q \in I^{[q]}$ in $R$.
 This implies that $\frac{z}{1} \in I^*R[U^{-1}]$ and the proof is complete.
\enddemoo

\demo{Proof of Theorem}
By the first lemma,
 it suffices to consider the case where $R$ is a domain.
The proof is complete by the second lemma.
\enddemoo

\proclaim{Corollary}
If $R$ has the property that modulo each minimal prime, the ring $R/P$ has an integral
extension that is F-regular, then tight closure commutes with localization in $R$.
\endproclaim

\demo{Proof}
By defintion, a ring is F-regular if all ideals are tightly closed both in the ring
and in all its localizations. Thus tight closure trivially commutes with localization
in an F-regular ring, and the corollary follows immediately from the Theorem.
\enddemoo

\proclaim{Corollary} Let $R$ be a ring of the form $k[x_1, \dots, x_n]/J$ where 
$J$ is generated by binomials in the $x_i's$ (a binomial is any polynomial of the 
form $\lambda x_1^{a_1}\dots x_n^{a_n} - \mu x_1^{b_1}\dots x_n^{b_n},$ where $\lambda$
and $\mu$ are elements of $k$. In particular, every monomial is a binomial in this sense.) 
Then tight closure commutes with localization in $R$.
\endproclaim

\demo{Proof}
There is no loss of generality in assuming the field is algebraically closed,
by an argument similar to the one used to prove Lemma 2 (we can apply  Lemma 2 directly
if $\bar k \otimes_k R$ is a domain where $\bar k$ is the algebraic closure of $k$).
In this case, the minimal primes of $J$ are themselves binomial ideals \cite{ES},
so by Lemma 1, we may assume that $R$ is a domain. 
Now, a ring which is the quotient of a polynomial ring by a prime binomial ideal 
is isomorphic to a subring of a polynomial ring $k[t_1, t_2, \dots, t_d]$
generated by monomials in the $t_i$'s, so we assume $R$ is such a monomial
algebra \cite{ES, 2.3}. The normalization $\tilde R$ of such a monomial algebra is again 
a monomial algebra.
On the other hand, a normal monomial algebra is a direct summand of the polynomial
overrring $k[t_1, \dots, t_d]$ (see, {\it e.g.\/} \cite{BH, Exercise 6.1.10}).
Since a direct summand of a regular ring (or even F-regular ring) is F-regular,
we see that $\tilde R$ is F-regular, so by the Corollary above, tight closure
commutes with localization in $R$.
\enddemoo

Special cases of this corollary, including the case where $J$ is generated by monomials, had
been observed using different methods by Will Traves, Moty Katzman, Keith Pardue and 
Karen Chandler, and Irena Swanson and myself. Katzman's work \cite {K} 
includes an interesting
approach to the localization problem using Gr\"obner basis
  in the
special case of binomial ideals. Aldo Conca \cite{C}, and Karen Chandler and Keith Pardue,
also made progress towards  this Gr\"obner 
basis approach in the binomial ideal case.  In \cite{SS}, a different
approach to the localization problem is considered, in which the growth of the
primary components of Frobenius powers is analyzed; the problem is solved for monomial
ideals in monomial rings. In \cite {T}, the monomial  case is solved using differential 
operators. The
thesis of Doug McCulloch \cite{Mc} does not directly address the localization  problem, but
contains a systemic study of tight closure in binomial rings.

\enddocument